\documentclass[letterpaper, 10 pt , conference]{ieeeconf}

\IEEEoverridecommandlockouts                              % This command is only needed if 
\usepackage{graphicx} % for pdf, bitmapped graphics files
\usepackage{epsfig} % for postscript graphics files
\usepackage{cite} % allows to put citations like [1]-[5]
\usepackage{amsmath} % assumes amsmath package installed
\usepackage{amssymb}  % assumes amsmath package installed
\usepackage{mathrsfs}
\usepackage{color}
\usepackage{nicefrac}
\usepackage[linesnumbered, ruled, vlined]{algorithm2e}
\usepackage{mathtools}
\mathtoolsset{showonlyrefs=true}
\usepackage{amsthm}
\graphicspath{{.}}

\theoremstyle{plain}

\usepackage{etoolbox}
\usepackage{float}
\usepackage{tabularx}

\def\th@definition{
	\thm@headfont{\itshape} % Heading font is italic
	\thm@notefont{} % Note is same as heading
}
\makeatother

\newtheorem{dfn}{Definition}

\newtheorem{lemma}{Lemma}

\newtheorem{problem}{Problem}
\newtheorem{theorem}{Theorem}
\newtheorem{remark}{Remark}
\newtheorem{assumption}{Assumption}
\newtheorem{corollary}{Corollary}
\newtheorem{prop}{Proposition}

\newcommand{\R}{\mathbb{R}}
\newcommand{\N}{\mathbb{N}}

\newcommand{\numticks}{\kappa_t}

\makeatletter
\def\endthebibliography{%
	\def\@noitemerr{\@latex@warning{Empty `thebibliography' environment}}%
	\endlist
}
\makeatother

\makeatletter
\patchcmd{\@makecaption}
{\scshape}
{}
{}
{}
\makeatother

\title{\LARGE \bf
	Technical Report: A Totally Asynchronous Algorithm for Tracking Solutions to Time-Varying Convex Optimization Problems
}

\author{Gabriel Behrendt and Matthew Hale$^{\ast}$% <-this % stops a space
	\thanks{$^{\ast}$Gabriel Behrendt and Matthew Hale are with the Department of  Mechanical and Aerospace Engineering at the University of Florida, Gainesville, FL USA. Emails: \texttt{\{gbehrendt,matthewhale\}@ufl.edu.} 
	This work was supported by AFOSR under grant FA9550-19-1-0169
	and by ONR under grants N00014-19-1-2543 and
	N00014-21-1-2495. 
	}
}

\begin{document}
	\maketitle
	\begin{abstract}
This paper presents a decentralized algorithm for a team of agents
to track time-varying fixed points that are the solutions to 
time-varying convex optimization problems.
The algorithm is first-order, and
it allows for \emph{total} asynchrony
in the communications and computations of all agents, i.e., 
all such operations can occur with arbitrary timing
and arbitrary (finite) delays. Convergence
rates are computed in terms of the communications and computations that
agents execute, without specifying when they must occur.
These rates apply to convergence to the minimum of each individual objective function,
as well as agents' long-run behavior as their objective functions change. 
Then, to improve the usage of limited communication and computation
resources, we optimize the timing of agents' operations relative to changes
in their objective functions to 
minimize total fixed point tracking error over time. Simulation results
are presented to illustrate these developments in practice
and empirically assess their robustness to uncertainties 
in agents' update laws. 
\end{abstract}
	\section{Introduction}

%P1: large scale Convex programs and time-varying problems
%P2: Asynchronous algorithms are better than averaging
%P3: Block-based algorithms are good
%P4: Decentralized Time-varying optimization
%P5: "In this paper..." we optimize time-varying objectives totally asynchronously
%P6: Our contributions
%P7: This paper is organized as follows

Time-varying convex optimization problems arise in 
%Time-varying convex optimization has gained interest in the fields of 
machine learning, signal processing, robotics, power systems, and others~\cite{fosson2018,Jakubiec2013,fosson2021,arslan2016,hauswirth2018}. 
%\mh{Are you certain these all have time-varying optimization in them?}
%\gb{[1] online convex optimization in dynamic environments in ML, [2] maximizing a Posteriori Probability of a time-varying signal, [3] sparse time-varying optimization to track time-varying signal, [4] solves the time-varying optimization problem~$\text{argmin} \ \frac{1}{2} \Vert x - x_{d}(t) \Vert^{2}$ and projecting it onto "collision-free" local workspace that is time-varying in robotics, [5] Time-varying Projected Dynamical Systems with Applications to Feedback Optimization of Power Systems}
Such problems can model, for example, time-varying demands in power systems~\cite{tang2017}, and uncertain, dynamic environments for real-time control of autonomous vehicles~\cite{zheng2020}. 
%Such problems arise in when operating environments are highly dynamic such as signal detection [CITATION], formation control [CITATION], state estimation [CITATION], and online learning [CITATION]. 
%Due to the time variation of these problems, it is necessary to solve them with an online algorithm that can handle large amounts of data. 
Often, these problems are solved with networks of agents, e.g., the robots that are given a time-varying task specification
or a network of processors across which computations have been parallelized. 
Such instances then require a multi-agent algorithm to solve time-varying problems.

Many multi-agent systems face asynchrony in agents' communications and computations. For example, asynchronous communications can result from
environmental hazards that impede communications and/or interference among communication signals.
Asynchronous computations onboard agents can result from clock mismatches and different computational hardware.
As a result, agents may generate and share information with unpredictable timing. 
Asynchrony has been studied for static optimization problems
in numerous ways, e.g., \cite{ye2018,nedic2011,iutzeler2013,duchi2011,tsianos2012}. 
In particular, \cite{Bertsekas89} establishes conditions for the convergence of static
optimization under \emph{total} asynchrony, namely communications and computations
with arbitrary delays. To the best of our knowledge, total asynchrony
remains unaddressed in time-varying optimization, and
this paper fills this gap. 

The algorithm we present uses totally asynchronous gradient descent to 
track time-varying fixed points that are the solutions to time-varying
convex optimization problems. It is block-based, in the sense that
each agent updates and communicates only a small subset
of the decision variables in a problem. Due to total asynchrony,
these operations can occur with arbitrary timing, which precludes
the development of convergence rates that depend on delay bound
parameters \cite{chang2016,mansoori2017,tseng1990,tsitsiklis86}.  
In this work, convergence rates are instead derived
that depend upon how agents' communications and computations
are interleaved over time. 
These results enable the quantification of the quality of agents'
fixed point tracking in terms of what they execute,
without specifying when operations must occur.

While it may be difficult to control the exact timing of all agents'
operations, it is easier for individual agents to control
the sequence of their own local operations (without needing
to control their exact timing). 
To leverage this capability, the aforementioned convergence
rates are used to optimize the sequence of 
agents' operations.
This is done by optimizing (in a way made
precise in Section~\ref{sec:Cycles}) the order
of communications and computations, both relative to each other
and relative to the timing of changes in objective functions.
%These results account for cases in which limited resources,
%such as finite battery power, place hard limits on the
%number of computations and communications that agents can execute. 
Specifically, a constrained optimization problem is formulated and solved in closed
form to minimize agents' convergence rates (through their
dependence on agents' operations) under constraints
on communicating and computing. 

\noindent \textbf{Summary of Contributions:}
In summary, the contributions of this paper are:
\begin{itemize}
	\item We show that the totally asynchronous gradient descent algorithm can track
	%time-varying fixed points that are the 
	solutions to time-varying convex optimization problems (Section~\ref{sec:results}).
	%and we quantify tracking error in terms of agents' communications and computations. 
	%and we implement it over a network of agents that communicate and compute updates independently. 
	\item We derive tracking bounds of the proposed algorithm in both the asymptotic and finite-time cases. These bounds 
	quantify the dependence of agents' performance upon their communications and computations, as well as problem
    parameters (Theorem \ref{thm:thm2} and Corollary \ref{cor:cor1}). 	
	%except for the communication between agents. This characteristic of tunable communication leads to our final contribution. \mh{Do the same here.} \gb{I tried to transition to the final contribution, don't know if it helps}
	\item For problems in which agents can execute a limited number of computations and communications,	
	we optimize the sequence of agents' operations to achieve desired tracking error.
	Namely, a constrained optimization problem is solved in closed form to determine the optimal
	sequence of communications and computations to enforce a desired bound on 
    fixed point tracking error (Theorems~\ref{thm:thm3} and \ref{thm:thm4} and Corollary~\ref{cor:cor2}). 
    \item We empirically demonstrate
    the robustness to noise of the totally asynchronous gradient descent algorithm in
    a time-varying problem (Section~\ref{sec:Simulation}). 
\end{itemize}

\noindent \textbf{Related Work} 
There is a large corpus of research on multi-agent optimization, and we cite only
a few representative examples. 
Asynchronous block-based algorithms include seminal work in~\cite{tsitsiklis1986, Bertsekas89, bertsekas1983}
and recent developments in \cite{hendrickson2020, ubl2020,yazdani2022,assran20,feyzmahdavian21},
which all consider static optimization problems. We differ by studying time-varying
problems. 
%This paper also differs from consensus optimization algorithms \cite{zhu2012,nedic2009,lobel2011},
%in which all agents update all decision variables and repeatedly average
%their iterates over time. We also do not assume that
%the agents' objective function has any particular form, e.g., being
%the sum of local objective functions, and instead we allow it to be arbitrary and arbitrarily coupled,
%subject to mild conditions.

%This paper is also closely related to a growing body of research on \emph{feedback optimization}, which
%generally considers time-varying optimization problems, often with some form
%of state measurements of a dynamical system entering computations. 
%Existing work most closely related to the present work pertains to the development of prediction-correction algorithms for problems with time-varying objective functions in both the centralized \cite{simonetto2016,simonetto2017} and distributed cases \cite{bernstein2018, simonetto2015}. 
%Specifically in \cite{simonetto2015}, they consider an objective function that is the sum of locally available functions at each node and uses higher-order information of their problem to perform the prediction step of their algorithm. In the present work we do not consider a prediction step in our algorithm because we consider an objective function that is globally available to every agent, not the sum of individual objective functions. Also, we only use first-order information in our algorithm because higher-order information can be hard to compute or unavailable.
This paper is also closely related to a growing body of research on \emph{feedback optimization}, which
often includes time-varying optimization problems. 
Relevant work develops prediction-correction 
algorithms for problems with time-varying objective functions in both the centralized \cite{simonetto2016,simonetto2017} 
and distributed cases \cite{simonetto2015,bernstein2018}. 
Specifically, \cite{simonetto2015} considers an objective function that is the sum of locally available functions
at each node and uses higher-order information to perform the prediction step.
The present work considers objective functions of an arbitrary form and requires only first-order information. 
Moreover, the objective functions we consider 
change discretely in time
%we make no assumptions about how they depend on time. 
%In particular, they 
%need not be continuous in time, 
and they need not be sampled from an objective
function that varies continuously in time. 
Changes in the objective function thus
cannot be predicted reliably, and we do not use a prediction step. 
In \cite[Section 1]{simonetto2015}, it is noted
that correction-only algorithms (which we consider)
cannot effectively mitigate tracking error of the optimizer 
if the optimization problem and algorithm do not have sufficient timescale separation. 
In this work, some timescale separation is indeed required, 
%We 
%discuss the timescales we consider in Section~\ref{sec:Update}, and indeed some timescale
%separation is required 
in the sense that objective functions must not vary in time
faster than agents can perform computations and communications 
(Section~\ref{sec:Update} provides a detailed discussion). 
We note that the need for timescale separation
is motivated by feasibility rather than technical convenience: if objectives can change faster than
agents can compute and communicate, then agents will be unable to track time-varying
fixed points and we cannot hope to prove anything useful. 

In~\cite[Remark 3]{bernstein2018}, 
it is noted that distributed time-varying optimization by agents with different computational abilities is 
subject to ongoing research. To the best of our knowledge, this problem remains open, and we 
present a solution to it by allowing agents to both compute and communicate totally asynchronously.

The rest of the paper is organized as follows. Section~\ref{sec:problemFormulation} states the problems of interest. 
Section~\ref{sec:Update} states the algorithm for solving them. 
Then, Section~\ref{sec:results} develops convergence rates for agents' short-term and asymptotic behavior. 
Section~\ref{sec:Cycles} optimizes agents' communications and computations to improve tracking over time. 
After that, Section~\ref{sec:Simulation} provides simulation results, and
Section~\ref{sec:conclusions} concludes.

	\section{Problem Formulation}\label{sec:problemFormulation}
This section gives a formal problem statement
and the assumptions placed on it. 
Throughout, this paper
solves problems with a network of~$N$ agents indexed
over the set~$[N] := \{1, \ldots, N\}$. 
We consider problems of the following form over a time horizon~$\mathcal{T} = \{0, \ldots, T\}$ for $T \in \N$.

\begin{problem} \label{prob:main}
With~$N$ agents, over a time horizon~$\mathcal{T}$,
\begin{equation}
	\begin{split}
		\underset{u}{\textnormal{minimize}} & \ f(u,t) \\
		\textnormal{subject to} & \ u\in U, 
	\end{split}
\end{equation}
where $U \subseteq \R^n$, $f:\R^{n}\times\N \rightarrow \R$, and $t\in \mathcal{T}$ is time. 
\end{problem}

We make the following assumptions about Problem \ref{prob:main}.

\begin{assumption} \label{as:c2}
For all~$t \in \mathcal{T}$,
	the function $f(\cdot, t)$ is twice continuously differentiable.
\end{assumption}

Assumption~\ref{as:c2} guarantees the existence and continuity of both the gradient and
 the Hessian of~$f(\cdot, t)$. 
 %that~$f$ has both a gradient and Hessian with respect to~$u$. 
 Many functions satisfy this assumption, such as polynomials of all orders.

\begin{assumption} \label{as:u}
	The constraint set $U$ can be decomposed via	
	\begin{equation}\label{eq:constraintSet}
		U=U_{1}\times U_{2}\times\cdots\times U_{N},
	\end{equation}
	 where 
	 $U_i \subseteq \R^{n_i}$,~$n_i \in \mathbb{N}$, is non-empty, compact, and convex.
\end{assumption}

Assumption~\ref{as:u} permits constraints such as box constraints, and it implies
that~$U$ is also non-empty, compact, and convex.
It also ensures that each agent can update its block of decision
variables independently, which will be used to develop a distributed gradient projection
law for enforcing the set constraint~$u \in U$. Similarly, $u\in U$ can be decomposed via
 
\begin{equation}\label{eq:partition}
    u= \left(\begin{array}{c}
        u_{1}\\
        \vdots\\
        u_{N}
    \end{array}\right), 
\end{equation}
  where $u_{i}\in U_{i} \subseteq \mathbb{R}^{n_i}$. 
  We define 
  \begin{equation} \label{eq:pu}
  p_u = [n_1, n_2, \ldots, n_N] \in \mathbb{N}^N
  \end{equation}
  to be the vector of dimensions of the blocks
  into which~$u$ is partitioned. 
  From Assumptions~\ref{as:c2}-\ref{as:u}, we have that~${H(\cdot,t) \coloneqq \nabla^{2} f(\cdot,t)}$ is continuous and~$U$ is compact.
  Therefore,~${\nabla f}(\cdot, t)$ is Lipschitz continuous. 
  For each~$t \in \mathcal{T}$, 
  we use~$L_t$ to denote the bound on~${\|H(\cdot, t)\|_2}$ over~$U$,
  and this~$L_t$ is the Lipschitz constant of~${\nabla f}(\cdot, t)$ over~$U$. 

The next assumption uses strict block diagonal dominance of partitioned matrices. 
Suppose a matrix $A \in \R^{n \times n}$ is partitioned into blocks according to the partition 
vector ${p_{u} = [n_1, \ldots, n_N]}$,
%\gb{Can we just use~$p_{u}$ here instead of defining a new~$p$?}, 
%\mh{Yes, based on the changes we've made, it makes sense to make this all about~$p_u$ now. 
%Before we did have an arbitrary~$p$ in our definitions and then we specialized to~$p_u$, but now
%it should all be~$p_u$. I've done some of it here, but propagate this through
%the rest of the paper. \gb{done}
%}
where $\sum_{i=1}^{N} n_i = n$ (we allow $n_i \neq n_j$ for all $i \neq j$). 
Let~$\left[ A \right]_{p_u}$ denote the partitioned matrix, i.e.,
%With this notation in place we can write a partitioned matrix as follows
\begin{equation}\label{eq:dd_def}
	\left[ A \right]_{p_{u}}=
	\begin{bmatrix}
	    A_{11} & A_{12} & \cdots & A_{1N}\\
		A_{21} & A_{22} & \cdots & A_{2N}\\
		\vdots & \vdots & \ddots & \vdots\\
		A_{N1} & A_{N2} & \cdots & A_{NN}
	\end{bmatrix}.
\end{equation}
The $i^{th}j^{th}$ block, denoted~$A_{ij} \in \R^{n_i \times n_j}$, is the matrix block of the rows
of~$[A]_{p_{u}}$ 
with indices~$\sum_{k=1}^{i-1} n_k + 1$ through $\sum_{k=1}^{i} n_k$ and
the columns of~$[A]_{p_{u}}$ with indices~$\sum_{k=1}^{j-1} n_k + 1$ 
through~$\sum_{k=1}^{j} n_k$, i.e.,~$A_{11}$ is the first~$n_1$ columns of the first~$n_1$ rows,~$A_{12}$
is the next~$n_2$ columns of the first~$n_1$ rows, etc.

%\mh{Make sure we change~$[A]_p$ to~$[A]_{p_u}$.}
%\gb{done}

\begin{assumption} \label{as:diag}
    The matrix
	$H(u, t)$ can be partitioned with respect to~$p_u$ from Equation~\eqref{eq:pu}, 
	and $H_{ii}(u,t)$ is symmetric and positive definite for all $u\in U$ and each $t \in \mathcal{T}$, i.e.,	
	%\begin{equation}\label{eq:sym_pd}
	$H_{ii}(u,t) = H_{ii}(u,t)^{T} \succ 0$.
	%\end{equation}
	Moreover, for each~$t \in \mathcal{T}$,
	$H(u,t)$ is $\beta_t$-strictly block diagonally dominant, i.e., for some~$\beta_t > 0$
	\begin{equation}\label{eq:ddH}
		\left( \left\Vert H_{ii}(u,t)^{-1}\right\Vert _{2} \right)^{-1} \geq \underset{j \neq i}{\sum_{j=1}^{n}} \left\Vert H_{ij}(u,t)\right\Vert _{2}+\beta_{t}.
	\end{equation}		
\end{assumption}

%\mh{Assumption~3 basically restates Definition~1, which I think is good. But then we can probably cut Definition~1. Do we use it anywhere
%apart from Assumption~3? If not, we can cut Definition~1.}
%\gb{I cut Def 1.}
%\mh{Assumption~3 uses~$[T]$. 
%However, it isn't defined yet. One thing you could do is, in the first paragraph
%of this section (by the definition of~$[N]$), you could say ``Problems evolve over the time horizon~$[T] = \{1, \ldots, T\}$
%for some~$T \in \mathbb{N}$.''
%}
%\gb{done}

As noted in~\cite[Section 6.3.2]{Bertsekas89},
Assumption~\ref{as:diag} is sufficient to ensure that 
agents' iterates contract towards the minimizer of~$f(\cdot, t)$
(in an appropriate norm) for fixed~$t$. 
And, in a way made precise in~\cite[Section 6.3.1]{Bertsekas89},
a problem is unlikely to
be solvable by a totally asynchronous algorithm 
 if it fails to satisfy Assumption~\ref{as:diag}. 
Thus, we enforce it here. 
%Assumption~\ref{as:diag} can be satisfied by a number of function families 
%such as... \gb{WHAT ARE FUNCTIONS THAT SATISFY THIS?}. 
One consequence of Assumption~\ref{as:diag} is given next. 
%Given our interest in developing a totally asynchronous 
%algorithm that outperforms its synchronous counterpart, \emph{Assumption 4} allows us to do so with 
%convergence guarantees.

\begin{lemma} \label{lem:Gershgorin}
	Let Assumptions~\ref{as:c2}-\ref{as:diag} hold. Then~$f(\cdot, t)$ is 
	$\beta_{t}$-strongly convex and it has a unique minimizer, denoted~$\hat{u}^t$.
\end{lemma}
\emph{Proof:} 
	From~\cite[Lemma 4]{ubl2020}, it immediately follows that
	\begin{equation}
	    \lambda_{min} \big[H(u, t)\big] \geq \big(\Vert H_{ii}(u,t)^{-1} \Vert_{2} \big)^{-1} - \underset{j \neq i}{\sum_{j=1}^{n}} \left\Vert H_{ij}(u,t)\right\Vert_{2}
	\end{equation}	
	for all~$u \in U$ and~$t \in \mathcal{T}$.
	%\mh{There is still~$\mathbb{N}^{+}$ here...Get rid of it everywhere.}
	%\gb{done}
	Therefore by Assumption~\ref{as:diag},
	\begin{equation}
	    \lambda_{min}\big[H(u,t)\big] \geq \beta_{t}
	\end{equation}	
	for all~$u \in U$ and~$t \in \mathcal{T}$. This implies that~$f(\cdot, t)$ is~$\beta_t$-strongly convex~\cite[Section 9.1.2]{boyd2004},
	from which the uniqueness of the minimum~$\hat{u}^t \coloneqq \underset{u \in U}{\text{argmin}} \ f(u,t)$ follows. 
\hfill $\blacksquare$

Regarding changes over time, the following is assumed.
 
\begin{assumption} \label{as:uhat}
    For each $t \in \mathcal{T}$, there exists a constant ${\sigma_{t} \geq 0}$ such that 
	%\begin{equation}\label{eq:drift_bound}
		$\left\Vert \hat{u}^{t+1}-\hat{u}^{t}\right\Vert_{2} \leq \sigma_{t}$.
	%\end{equation}
\end{assumption}

%Furthermore, we define the set ~$[T] \coloneqq \{1,\ldots,T\}$.
Without Assumption~\ref{as:uhat}, a problem 
could change by arbitrarily large amounts, and
it would not be possible to track the solution of Problem~\ref{prob:main}. 

	\section{Asynchronous Update Law}\label{sec:Update}
    %\mh{I think we should introduce the idea of timescale separation here. Say that
    %changes in the optimization problem are indexed over~$t$ and we assume that
    %some amount of clock time elapses between switches in the problem.
    %Then say that that clock time enables agents to perform some number of operations,
    %which we index over the counter~$k$. 
    %Then introduce an assumption that, between~$t$ and~$t+1$ in the problem clock,
    %some number of ticks~$k_t \geq 1$ elapse for the agents' clock. 
    %Then say that agents may compute something new or communicate at each such tick,
    %though this is not guaranteed due to asynchrony. Talk about how this must be here for
    %feasibility: if agents don't have time to compute or communicate anything, then
    %they obviously can't track fixed points.
    %}
    %\gb{done}
    
    In this section we develop the totally asynchronous
    algorithm for tracking time-varying fixed points.
    First, we describe the relationship between
    the times at which the objective function
    changes and the times at which agents execute communications
    and computations.
    
    \subsection{Time and Timescale Separation}
    Problem~\ref{prob:main} evolves over time, indexed by the discrete time index~$t$, and, in
    defining agents' algorithm, we will index their operations over
    the discrete time index~$k$. The reason is that the timing of agents'
    computations and communications need not coincide with changes
    in their objective functions. For example, agents may execute
    several communications and computations when minimizing~$f(\cdot, t)$
    for each~$t \in \mathcal{T}$. If this were not done, 
    then the problem could evolve faster than agents compute and communicate, and
    tracking a time-varying fixed point would be infeasible.
    Thus, it is necessary to adopt a formalism that enables agents to perform
    multiple operations to minimize~$f(\cdot, t)$ for each~$t \in \mathcal{T}$.
    
    In particular, we require some
    timescale separation between the ticks of~$t$, the time index
    for changes in the objective function,
    and the ticks of~$k$, the index of agents' operations. For
    for each~$t \in \mathcal{T}$, 
    some number of ticks~$\numticks \geq 1$ of~$k$ elapse 
    before~$t$ increments to~$t+1$.
    Agents may compute and/or communicate something at each tick of~$k$,
    though this is not guaranteed due to asynchrony. 
    In other words, timescale separation allows for the possibility that
    agents perform multiple communications and computations
    for each~$t$, though we by no means assume that they do so.
    As will be shown below, our analysis allows for (and quantifies
    the effects of) scenarios in which agents perform no operations at all.
        Of course, timescale separation
    is not new to this work, and has been used, e.g., in optimization and
    model predictive control in~\cite{tan2016,weiss2014,hauswirth2021}. 
    %\mh{Cite MPC works,
    %and see if it's in that other distributed asynchronous
    %time-varying fixed point tracking paper (it should be in there).
    %It may be implicitly in there rather than explicitly like we're doing,
    %but let's briefly mention it here for sure. 
    %}
    %\gb{28-29 are MPC, 30 is feedback optimization with timescale separation. Other distributed TV tracking paper says the following about timescale "Remark 3: The considered setting does not capture the
    %case where agents exhibit different computational capabilities. The analysis of a more general case where asynchronous
    %updates are due to both communication constraints and
    %different computational times is the subject of ongoing work"}
    
    \subsection{Formal Algorithm Statement}
	We consider a block-based gradient algorithm with totally asynchronous computations and communications to solve Problem~\ref{prob:main}	
	over a network of~$N$ agents. Each agent updates only a subset of the entire decision vector and communicates 
	the values of their sub-vector with the rest of the network over time. 
	Asynchrony implies that agents receive different information at different times, and thus we expect
	them to have differing values for decision variables onboard. 
	At any time $k$, agent $i \in \left[ N \right]$ has a local copy of its decision vector, denoted $u^{i}(k) \in \R^{n}$. Due to 
	computation and communication delays, we allow $u^{i}(k) \neq u^{j}(k)$ for~$j \neq i$. 
	Within the vector~$u^i(k)$, agent~$i$ computes updates only to~$u_{i}^{i}(k) \in \R^{n_i}$.  
	%where $n_i \in \mathbb{N}$ is the length of agent $i$'s sub-vector. 
	At any time $k$, agent $i$ 	has (possibly old) values for agent $j$'s sub-vector of the decision variables, 
	denoted~$u_{j}^{i}(k) \in \R^{n_j}$. 	
	These are updated only by agent~$j$ and then shared with the rest of the network. 

    %Now we discuss the communication and update law of the agents. 
    At time $k$, 
    if agent~$i$ computes an update to~$u^i_i(k)$, 
    its computations are in terms of
    the decision vector $u^{i}(k)$. The entries of $u_j^{i}(k)$ 
    for~$j \neq i$
    are obtained by communications 
    from agent~$j$, which may be subject to delays. 
    Therefore, agent~$i$ may (and often will) 
    compute updates to~$u^i_i$ based on outdated information from
    other agents. 
    
    Let $K^{i}$ be 
    the set of times at which agent~$i$ computes an update to~$u^i_i$.
    Similarly, let~$C^{i}_j$ 
    be the set of times at which agent~$j$ receives a transmission
    of~$u^i_i$ from agent~$i$; due to communication delays, these transmissions can be received at any time
    after they are sent, and they can be received at different times by different agents. 
    We emphasize that the sets~$K^i$ and~$C^i_j$ are only defined
    to simplify discussion and analysis. Agents do not know and do not need
    to know~$K^i$ and~$C^i_j$. 
    %\mh{I've changed~$C^i$ to~$C^i_j$, and it's now the set of times when agent~$j$ receives
    %something from agent~$i$. Before, it seemed like communications were synchronous
    %because there was only one~$C^i$. Propagate this change forward to each
    %other place that needs it. I've included this in Algorithm~1 already.}
    %\gb{done. Looks like this is the only spot it is used}
    We define~$\tau_{j}^{i}(k)$ to be the time at which agent~$j$ originally computed the value of~$u^i_j$ that
    agent~$i$ has onboard at time~$k$; agents~$i$ and~$j$ also do not need to know~$\tau^i_j$.    
    Formally, 
    \begin{equation}
    u^i_j(k) = u^j_j\big(\tau^i_j(k)\big). 
    \end{equation}
    %because there can be a delay between computation and communication of sub-vectors. 
    Then~$k-\tau_{j}^{i}(k)$ is the length of communication delay, and we
    emphasize that \emph{we do not assume any bound on agents' delays}. 
    Instead, only the following is enforced.
        
    \begin{assumption} \label{as:tau}
        For all~$i \in [N]$, the set~$K^i$ is infinite. Let~$\{k_d\}_{d \in \N}$ be a non-decreasing sequence
        in~$K^i$ tending to infinity. Then for all~$i \in [N]$ and~$j \in [N] \backslash \{i\}$, we have
        \begin{equation}
            \lim_{d \to \infty} \tau^i_j(k_d) = \infty. 
        \end{equation}        
        And messages are received in the order they were sent.
    \end{assumption}
    
    Assumption~\ref{as:tau} ensures that (i) no agent ever stops computing updates (because~$K^i$ is infinite),
    (ii) no agent ever stops sharing its updates with other agents (because~$\tau^i_j(k_d)$ tends to infinity),
    and (iii) arbitrarily old information cannot be received arbitrarily late into the run of an algorithm.
    This assumption is one of feasibility: there is no tracking of fixed points (and nothing to prove)
    if any agent ever permanently stops computing and/or communicating, or if arbitrarily old
    information can be received at any time. 
    
    We propose to use a gradient-based update law for agent~$i$, which takes
    the form
    \begin{align}
    	u_{i}^{i}(k+1) &=
    	\begin{cases}
    		\Pi_{U_i}\!\left[u_{i}^{i}(k) \!-\! \gamma_{t}\left(\nabla_{u_{i}}f(u^{i}(k),t) \right)\right] & k\!\in\! K^{i}\\
    		u_{i}^{i}(k) & k \!\notin\! K^{i}
    	\end{cases} \\
		u_{j}^{i}(k+1) &=
		\begin{cases}
			u_{j}^{j}(\tau_{j}^{i}(k)) & k \in C^j_i\\
			u_{j}^{i}(k) & k \notin C^j_i
		\end{cases},
	\end{align}

	\noindent where $\gamma_{t} > 0$ is the stepsize used when minimizing~$f(\cdot, t)$. 
	Bounds on~$\gamma_{t}$ are established later in Proposition~\ref{prop:TI_convergence}. 
	Algorithm~\ref{alg:myAlg} provides pseudocode for agents' updates; for notational convenience, we set~$\kappa_{-1} = 0$.
	%\mh{for notational convenience, we set~$\kappa_0 = 0$.} \gb{done}

%\mh{Should lines 4-9 be inside another ``for'' loop? I think it should be over~$j=0:N$. We use the index~$j$ in those lines,
%but it's not actually defined inside the algorithm.}
%\gb{yup. I changed it.}
%\mh{Check where these indices all start. I think~$i$ and~$j$ should start at~$1$, not~$0$.}
%\gb{yes, i agree. I updated that in the algorithm}

\begin{algorithm} \label{alg:myAlg}
	\SetAlgoLined
	\KwIn{$u_{j}^{i}(0)$ for all~$i,j \in [N]$,~$T \in \mathbb{N}$}
	
	$k_{loop,-1} = -1$
	
	\For{t=0:T}
	{
	    $k_{loop,t} = k_{loop,t-1} + \kappa_{t-1}$
	    
		\For{$k=k_{loop,t}+1:k_{loop,t}+\kappa_{t}$}
		{
			\For{i=1:N}
			{
			    \For{j=1:N}
			    {
    				\If{$k \in C^{j}_i$}{    			
    					~$u_{j}^{i}(k+1) = u_{j}^{j}(\tau^i_j(k))$ %\mh{Don't we need a~$\tau^i_j(k)$ in here? This makes comms. seem instant.}
    					%\gb{done}
    				}
    				\Else{
    				~$u_{j}^{i}(k+1) = u_{j}^{i}(k)$
    				 %\mh{Else it stays constant. Write this in math.} \gb{done}
    				}
    			}
				\If{k $\in K^{i}$ }{
					\!\!\!${u_{i}^{i}(k\!+\!1) \!=\!
					    \Pi_{U_i}\!\left[u_{i}^{i}(k) \!-\! \gamma_{t}\left(\nabla_{u_{i}}f(u^{i}(k),t) \right)\right]}$
				}
				\Else{
				~$u_{i}^{i}(k+1) \!=\! u_{i}^{i}(k)$ %\mh{Else it stays constant. Write this in math.}\gb{done}
				}
			}			
		}	
	}
	\caption{Totally Asynchronous Algorithm for Tracking Time-Varying Fixed Points}
\end{algorithm}

The remainder of the paper focuses on analyzing and optimizing
the execution of Algorithm~\ref{alg:myAlg}.

%	We note that since \eqref{eq:Opt} varies discretely over $t$ that we do require the agents to complete at least one communication cycle every $t$. This is a mild assumption because stated simply we require that the network makes some progress toward the minimizer $\hat{u}^t$ for all $t$.
	
\section{Convergence Results}\label{sec:results}
This section establishes error bounds
for Algorithm~\ref{alg:myAlg} 
by first deriving convergence rates
for time-invariant problems. Then, those results are extended to the time-varying case
to bound error when tracking time-varying fixed points. 
%This is possible by observing that solving a time-varying problem in discrete time is 
%essentially tracking the solution of a series of time-invariant problems. 
First, we 
define block-maximum norms. 

\begin{dfn}[Block-maximum norm] \label{def:max_norm}	
	Let $u \in \R^{n}$ be partitioned according to~$p_{u}$ as in Equation~\eqref{eq:pu}. Then define
	\begin{equation}
		\left\Vert u \right\Vert _{2,\infty}	\coloneqq \max_{i \in [N]} \left\Vert u_i \right\Vert_{2}.
	\end{equation}

\end{dfn}

This definition allows for the quantification of the worst-performing agent in the network, which will be considered in the convergence analysis.

For a time-invariant objective function, 
it has been shown in~\cite[Section 6.2, Prop. 2.1]{Bertsekas89} 
that if one can construct a 
sequence of sets~$\left\{U(s)\right\}_{s \in \N}$ that satisfy the properties of 
Definition~\ref{def:sets} below, then 
Algorithm~\ref{alg:myAlg} asymptotically converges to a minimizer of
that objective function. 
We will leverage that analysis for each fixed~$t$ and each fixed~$f(\cdot, t)$ individually. 
Then we will use these analyses to bound tracking error when~$f$ evolves over time. 
%We first focus on minimizing~$f(\cdot, 0)$ \gb{changed this from~$f(\cdot, t)$ to $f(\cdot, 0)$}, then consider
%time-varying objectives. 
For any fixed~$t$, we have the following. 

\begin{dfn} \label{def:sets}
    Fix~$t \in \mathcal{T}$ and fix~$f(\cdot, t)$. 
    A collection of sets~$\{U(s)\}_{s \in \N}$ is called \emph{admissible} if it satisfies:
    \begin{enumerate}
		\item For all~$s \in \N$, $U(s)$ can be decomposed as $U(s) =U_{1}(s) \times \cdots \times U_{N}(s)$,
		where~$U_i(s) \subseteq U_i$.
		\item $U(s+1)\subset U(s)$ for all~$s \in \N$
		\item $\underset{s\rightarrow\infty}{\lim} \ U(s)=\left\{ \hat{u}^{t}\right\}$
		\item For all $y \in U(s)$ and $i\in \left[ N \right]$, the update $z_{i}=\Pi_{U_{i}}\left[y_{i}-\gamma\nabla_{y_{i}}f(y,t)\right]$ gives $z_{i}\in U_{i}(s+1)$.
	\end{enumerate}
\end{dfn}

Definition~\ref{def:sets}.1 allows for the agents to independently update their local sub-vectors asynchronously
without violating set constraints. 
Definitions~\ref{def:sets}.2 and \ref{def:sets}.3 guarantee that the collection of sets $\left\{ U(s)\right\}_{s\in\N}$ is nested and converges to a singleton containing $\hat{u}^{t}$. Definition~\ref{def:sets}.4 guarantees that each update of an agent's sub-vector 
makes progress toward $\hat{u}^{t}$.

In the analysis of Proposition~\ref{prop:TI_convergence} below, we fix a~${t \in \mathcal{T}}$ 
in order to only analyze the convergence of Algorithm~\ref{alg:myAlg} 
for~$f(\cdot, t)$. 
After proving convergence to the optimizer of a fixed~$f(\cdot, t)$, 
we will extend these results to the time-varying case.
We introduce the symbol~$\eta_t\in \N$ for this analysis, and it denotes the 
total number of ticks of~$k$ that elapse before~$t$ increments to~$t+1$,
i.e.,~$\eta_t = \sum_{\ell=0}^{t} \kappa_t$. 

%\gb{we use $\eta_{t-1}$ here in Prop 1 but don't define it until theorem 2}
\begin{prop}\label{prop:TI_convergence}
    Let Assumptions \ref{as:c2}-\ref{as:uhat} hold. Fix~$t \in \mathcal{T}$ and fix~$f(\cdot, t)$. 
    Consider the collection~$\{U(s)\}_{s \in \N}$ 
	with
	\begin{equation}\label{eq:proposed_sets}
		U(s) \!=\! \left\{ y \in \R^{n}:\left\Vert y-\hat{u}^{t}\right\Vert _{2,\infty}\leq q_{t}^{s} \max\limits_{i \in [N]} \left\Vert u^{i}(\eta_{t-1}) - \hat{u}^{t} \right\Vert_{2,\infty}\!\right\},
	\end{equation}
	%\gb{do we need a~$u^{i}(\eta_{t})$ instead of~$u^{i}(0)$ since this is for any fixed~$t$???}
	%\mh{Good catch! I think it should be~$u^i(\eta_{t-1})$, since that's the last~$u^i$ generated before~$f(\cdot, t-1)$ changes to~$f(\cdot, t)$, and
	%that's the value of~$u^i$ that agent~$i$ has immediately after the switch (but before it computes anything). Let me know what you think.}
	%\gb{done}
	%\mh{Make sure the use of~$\eta_{t-1}$ shows up in the proof too (in the appendix).}
	where the stepsize~$\gamma_t > 0$ satisfies
	\begin{equation}
	\gamma_{t}<\frac{1}{\max\limits_{u \in U} \max\limits_{i \in [N]} \sum\limits_{\substack{j=1  \\ j \neq i}}^{n} \left\Vert H_{ij}(u,t)\right\Vert_{2}}, 
	\end{equation}
	and
	$q_{t} \coloneqq \max \{  \left| 1- \gamma_{t} \beta_{t} \right|, \left| 1-\gamma_{t} L_{t} \right| \}$,
	which satisfies~$q_t \in (0, 1)$. 
	This collection is admissible in the sense of Definition~\ref{def:sets}. 
\end{prop}
\emph{Proof}: See Appendix \ref{app:TI_convergence}. \hfill $\blacksquare$

This construction of sets will be used to derive a convergence rate in terms
of the number of communication cycles agents complete, defined next. 

\begin{dfn}[{\cite[Definition 4]{hale17}}]	
    A communication cycle has occurred after 
    (a)~for all~$i$, agent~$i$ has computed at least one update to~$x^i_i$, and then
    (b)~that updated value has been sent to and received by the other agents that need it. 
	%We claim that a communication cycle has occurred when every agent in the network has performed an update and communicated it with the rest of the network. 	
	The number of communication cycles started and finished when minimizing~$f(\cdot, t)$ 
	before switching to~$f(\cdot, t+1)$ is denoted by~${c(t)\in\mathbb{N}}$.	
\end{dfn}

Note that some agents may compute and share many updates within one cycle because a cycle is only completed
after \emph{all} agents have computed and shared at least one update.

The following constants are used in the forthcoming analysis:
%\mh{What is~$u^{i,0}$? That's in the definition of~$D_0$.}
%\gb{typo. fixed it.}
\begin{align}
q_{max} &\coloneqq \max_{t \in \mathcal{T}} q_{t} 
\label{eq:qmax}\\
D_{0} &\coloneqq \max_{i \in [N]} \left\Vert u^{i}(0)-\hat{u}^{0}\right\Vert_{2,\infty} \\
D_{t} &\coloneqq q_{t-1}^{c(t-1)}D_{t-1}+\sigma_{t-1}, \textnormal{ for all } t \geq 1 \label{eq:dtdef} \\
B &\coloneqq \max\left\{\max_{t \in \mathcal{T}} \sigma_{t},\ D_{0}\right\}.
\label{eq:Bdef}
\end{align}

We first analyze convergence for~$f(\cdot, 0)$ held fixed, then iterate
this analysis forward in time to account for time-varying problems. 
The convergence rate for~$f(\cdot, 0)$ is as follows.

\begin{theorem}[Time-Invariant Convergence Bound]\label{thm:thm1}
Let Assumptions~\ref{as:c2}-\ref{as:diag}
and~\ref{as:tau} hold. 
Then, for~$N$ agents executing Algorithm~\ref{alg:myAlg}, just before~$f$ switches from~$f(\cdot, 0)$ to~$f(\cdot, 1)$, 
for all~$i \in [N]$ we have
% For any~$t \in \mathcal{T}$ and~$\eta_t =\sum_{\ell=0}^{t} \kappa_{\ell}$, 

%\mh{What conditions is this under? What assumptions? What is the setup/context for this theorem?}
%\mh{I think this theorem needs to be
%for each~$t \in [T]$ and~$k = \sum_{\ell=0}^{t} \kappa_{\ell}$. 
%Then just say that this characterizes the tracking of~$\hat{u}^t$ just before~$f$ switches
%from~$f(\cdot, t)$ to~$f(\cdot, t+1)$.
%}
%\gb{done.}
    %\begin{equation}\label{eq:TI_tracking_bound}
	%	\left\Vert u^{i}(\eta_t)-\hat{u}^{t}\right\Vert _{2,\infty} \leq q_{t}^{c(t)}D_{t}.
	%\end{equation}
	\begin{equation}\label{eq:TI_tracking_bound}
		\left\Vert u^{i}(\eta_0)-\hat{u}^{0}\right\Vert _{2,\infty} \leq q_{0}^{c(0)}D_{0}.
	\end{equation}
\end{theorem}
%\mh{Careful with notation. The left-hand side had~$u^i(0)$ before, which was wrong.}
%\gb{Does this work at~$t=0$? $\eta_0 = \sum_{\ell=0}^0 \kappa_{\ell} = \kappa_0 = 0.$ So the bound at $t=0$ would be $\left\Vert u^{i}(0)-\hat{u}^{0}\right\Vert _{2,\infty} \leq q_{0}^{c(0)}D_{0}$ which I don't think is true}
%\mh{Yep, you're right. This happens because~$t$ starts at~$0$. But that's totally fine. I introduced~$\kappa_0 = 0$ because I thought~$t$ started at~$1$ before.
%So all that we would need to do is change~$\kappa_0 = 0$ at the end of Section III (it's in red)
%to~$\kappa_{-1} = 0$ to account for the fact
%that we use~$\kappa_{t-1}$ in Algorithm~1. I think that should fix it. What do you think?
%}
%\gb{yes I agree. And this works based on our definition of $\kappa_t$ in the timescale separation section.}
\emph{Proof:}
This follows from applying~\cite[Theorem 1]{ubl2020} to~$f(\cdot, 0)$ when agent~$i$ has the initial condition~$u^i(0)$. 
%\mh{We shouldn't prove this here. Just cite it from an earlier paper from the lab.}
%\gb{It seems like Ubl's proof in his paper "Totally Asynchronous Large-Scale..." %would be the closest to what I need here? Would you agree to cite that?}
\hfill $\blacksquare$
%\gb{need help relating theorem 2 to theorem 1 now per my slack message. I adjusted theorem 2 a little. Does it make sense based on the changes to theorem 1?}

%With a convergence rate for a fixed objective function,
%we can extend to the time-varying case. 
%\mh{What is~$u(k)$ without a superscript? Is it~$(u^1_1(k), \ldots, u^n_n(k))$? If so,
%say that because I don't believe that's been defined yet. 
%}
Using Theorem~\ref{thm:thm1},  the next result 
accounts for changes in the objective function and
bounds tracking error across finite time horizons
under such changes. 

\begin{theorem}[Finite-Time Tracking Bound]\label{thm:thm2}	
	%\mh{What are the assumptions/hypothesis of this result?} \gb{done}
	Let Assumptions~\ref{as:c2}-\ref{as:tau} hold. 
	Then, for any~$t \in \mathcal{T}$ and~$\eta_t = \sum_{\ell=0}^{t} \kappa_{\ell}$,
	the iterates produced by~$N$ agents executing 
	Algorithm~\ref{alg:myAlg} satisfy
	%\mh{On the LHS, is that~$u^i(t)$ or some other time?}
	%\gb{i believe it should be $u^i(k)$ so I changed it}
	\begin{equation}\label{eq:tracking_bound}
		\left\Vert u^{i}(\eta_t)-\hat{u}^{t}\right\Vert _{2,\infty} \leq D_{0} \prod_{\theta=0}^{t} q_{\theta}^{c(\theta)} + \sum_{p=1}^{t} \sigma_{p-1} \prod_{r=p}^{t}q_{r}^{c(r)}.
	\end{equation}
	%\mh{Be careful with~$k$. I think this holds for~$k = \sum_{\ell=0}^{t} \kappa_{\ell}$.}
	%\gb{I added this to the Theorem statement above}
\end{theorem}

\emph{Proof:} For all~$i \in [N]$, 
	%We proved the following in  due to our construction of sets in Equation \eqref{eq:proposed_sets}
	%\mh{I don't think we have that just yet. We need to explicitly state a convergence rate.}
	%\gb{I added Theorem 1 for the time-invariant case so I can reference it now.}
	Theorem~\ref{thm:thm1} gives
	\begin{equation}\label{eq:TI_tracking_bound}
		\left\Vert u^{i}(\eta_0)-\hat{u}^{0}\right\Vert _{2,\infty} \leq q_{0}^{c(0)}D_{0}.
	\end{equation}
	%for any~$t \in \mathcal{T}$ and~$\eta_t=\sum_{\ell=0}^{t} \kappa_{\ell}$. 
	%At~$t=0$ the bound in Theorem~\ref{thm:thm1} can be written as
	%\begin{equation}
	%    \left\Vert u^{i}(\eta_-1)-\hat{u}^{0}\right\Vert _{2,\infty}\leq q_{0}^{c(0)}D_{0}.
	%\end{equation}
	%\mh{For which values of~$k$?}
	%\mh{What is the time of~$u^i$ on the LHS?}
	%\gb{I think~$u^i(k)$}
	After the objective function changes from~$f(\cdot, 0)$ to~$f(\cdot, 1)$ (and before
	the next tick of~$k$), 
	the tracking error for agent~$i$ can be written as
	\begin{align}
	    \left\Vert u^{i}(\eta_0)-\hat{u}^{1}\right\Vert_{2,\infty} &\leq
	    \|u^i(\eta_0) - \hat{u}^0\|_{2, \infty} + \|\hat{u}^1 - \hat{u}^0\|_{2,\infty} \\
	    &\leq q^{c(0)}_0D_0 + \sigma_0,
	\end{align}
	which follows from the triangle inequality, Theorem~\ref{thm:thm1},
	Assumption~\ref{as:uhat}, and the fact that~$\|v\|_{2, \infty} \leq \|v\|_2$ for~$v \in \mathbb{R}^n$. 
	%\mh{For which values of~$k$?}
	%which can be re-written based on Equation \eqref{eq:dtdef} and the distance moved by the minimizer between switches in the objective function from~$f(\cdot,0)$ to~$f(\cdot,1)$ as
	Next,~$\kappa_1$ ticks of~$k$ will elapse before~$f(\cdot, 1)$ switches to~$f(\cdot, 2)$, during
	which time agents will complete some number of cycles, denoted~$c(1)$. 
	Then, just before~$f(\cdot, 1)$ switches to~$f(\cdot, 2)$, for all~$i \in [N]$ we find
	\begin{equation}
	    \left\Vert u^{i}(\eta_1)-\hat{u}^{1}\right\Vert _{2,\infty}\leq q_{1}^{c(1)}\left(q_{0}^{c(0)}D_{0}+\sigma_{0}\right).
	\end{equation}
%    \mh{For which values of~$k$?}	
	Applying this recursively to all~$t \in \mathcal{T}$ completes the proof.	
\hfill $\blacksquare$

%This bound is in terms of the maximum initial distance any agent starts from~$\hat{u}^{0}$, the distances travelled by the minimizer $\{\sigma_{t}\}_{t \in [T]}$, 
%the parameter $q_{t}$, and the number of communication cycles completed each time step $c(t)$. 
Increasing the number of communication cycles performed to minimize~$f(\cdot, t)$, namely~$c(t)$, will decrease the value of~$q_{t}^{c(t)}$, 
shrinking the error on the right-hand side of Equation~\eqref{eq:tracking_bound}. 
Therefore, the completion of more communication cycles will lead to better tracking of the minimizer~$\hat{u}^{t}$.
Indeed, the long-term tracking error depends on the value of~$c(t)$ for every previous~$t \in \mathcal{T}$,
which is intuitive because cycles count agents' total operations over time. 

The following result bounds tracking error as the length of time horizon~$T$ grows
arbitrarily large. It assumes that~$c(t) \geq 1$, i.e., that agents complete
at least one communicate cycle per objective function.
Without such an assumption, it is possible for agents to perform
no computations and communications, and no tracking would
be possible. 

\begin{corollary}[Asymptotic Tracking Bound] \label{cor:cor1}
Let Assumptions~\ref{as:c2}-\ref{as:tau} hold and suppose that N agents are executing Algorithm~\ref{alg:myAlg} with~$c(t) \geq 1$ for all~$t$. 
Then
\begin{equation}
    \underset{t \rightarrow \infty}{\lim}\left\Vert u^{i}(\eta_t)-\hat{u}^{t}\right\Vert _{2,\infty} \leq B \frac{q_{max}}{1-q_{max}}.
\end{equation}
\end{corollary}

\emph{Proof:}
%From Theorem~\ref{thm:thm2} we have
%\begin{equation}
%		\left\Vert u^{i}(\eta_t)-\hat{u}^{t}\right\Vert _{2,\infty} \leq D_{0} \prod_{\theta=0}^{t} q_{\theta}^{c(\theta)} + \sum_{p=1}^{t} \sigma_{p-1} \prod_{r=p}^{t}q_{r}^{c(r)}.
%\end{equation}
Equations~\eqref{eq:qmax} and~\eqref{eq:Bdef} imply that
\begin{align}
    B &\geq D_{0}, \quad B \geq \sigma_{t} \textnormal{ for all } t \in \mathcal{T}, \quad\textnormal{ and } \\
    q_{max} &\geq q_{t}  \textnormal{ for all }  t \in \mathcal{T}.
\end{align}
Then the result in Theorem~\ref{thm:thm2} can be upper bounded via
\begin{flalign}
    \left\Vert u^{i}(\eta_t)-\hat{u}^{t}\right\Vert _{2,\infty} &\leq B \prod_{\theta=0}^{t} q_{max}^{c(\theta)} + \sum_{p=1}^{t} B \prod_{r=p}^{t}q_{max}^{c(r)} \\
    &= B \left( \prod_{\theta=0}^{t} q_{max}^{c(\theta)} + \sum_{p=1}^{t} \prod_{r=p}^{t}q_{max}^{c(r)} \right). 
\end{flalign}
With~$c(t) \geq 1$ for all~$t \in \mathcal{T}$, we have 
\begin{flalign}
    \left\Vert u^{i}(\eta_t)-\hat{u}^{t}\right\Vert _{2,\infty} &\leq B \left( \prod_{\theta=0}^{t} q_{max} + \sum_{p=1}^{t} \prod_{r=p}^{t}q_{max} \right) \\
    & \leq B \left( \prod_{\theta=0}^{t} q_{max} + \sum_{p=1}^{t} q_{max}^p \right). 
\end{flalign}
The results follows by taking the limit as~$t \rightarrow \infty$ and using the fact~${\sum_{k=1}^{\infty} r^k = \frac{r}{1-r}}$ for~$\vert r \vert < 1$. 
\hfill $\blacksquare$

\begin{remark}
    A similar result to Corollary~\ref{cor:cor1} can be found in~\cite[Theorem 2]{bernstein2018}, which
    considers a synchronous model of computation. 
    Corollary~\ref{cor:cor1} thus shows that, under mild conditions, one can achieve
    similar tracking performance with a computational model that permits total
    asynchrony in computations and communications among agents.     
    Thus, total asynchrony can be permitted without sacrificing performance. 
\end{remark}

Resource-constrained agents may only be able to complete a limited number of total communications
and computations, e.g., due to limited fuel or battery power onboard a mobile robot.
The long-run dependence upon all values of~$c(t)$ then leads to the question
of how communications and computations should be interleaved over time
to tune the values of each~$c(t)$ to lead to optimal tracking overall.
This is the subject of the next section. 

%If agents are unable to complete any communication cycles between switches in~$f(\cdot, t)$, then~$c(t) = 0$.
%\mh{Say more. Is the point really that~$c(t)$ bigger means~$q_t^{c(t)}$ is smaller? If so,
%say that explicitly. We want to make sure we're telling the reader everything about our thought process,
%even if it seems obvious to us at this point.}
%\gb{done}
	
%\begin{corollary} [No Communication Tracking Bound] \label{cor:noComms}	
%	\begin{equation}
%		\left\Vert u^{i}(k)-\hat{u}^{t}\right\Vert _{2,\infty} \leq D_{0} + \sum_{p=1}^{t} \sigma_{p-1}
%	\end{equation}	
%\end{corollary}
%
%We note that the upper bound here is simply the sum of
%the initial distance to the minimizer with the sum of all distances
%it moves over time. In short, in the absence
%of any completed communication cycles, the distance to the optimizer 
%increases without bound. The attainment of any meaningful levels
%of performance therefore requires the completion
%of communication cycles over time. Noting that these cycles
%count agents' computations and communications interleaved in a
%certain way, we observe that the asynchronous execution
%of agents' operations will result in completed cycles over time,
%even without any deliberate interleaving of or coordination
%between these operations.

%For cases in which agents can choose when and how to
%interleave computations and communications, the next
%section optimizes their timing to ensure that
%agents' operations are as effective as possible
%in reducing tracking error. 

	\section{Cycle Optimization}\label{sec:Cycles}
%The previous section describes the tracking error bounds in terms of communication cycles and it was mentioned that these bounds can be made tighter by increasing the number of communication cycles completed by the network each time step. 
In this section, we determine when  agents  should compute and communicate in  order  to  minimize  the error bound
in Theorem~\ref{thm:thm2}. 
In particular, we optimize the number of communication cycles~$c(t)$ that agents should
perform when minimizing~$f(\cdot, t)$. This is done under constraints on the number
of communications and computations that agents can perform, which can represent
limits on total energy and/or bandwidth. 
%Determining the number of communication cycles required for certain performance specifications 
%is important because in most time-varying problems we do not have control over how the problem 
%evolves over time. However, we can determine how often agents compute updates and communicate 
%based on onboard processor specifications and environmental factors.
This is done for two cases: (i) when~$c(t) \equiv c$, i.e., when the number of cycles is
the same for all~$t$, and (ii) when the number of cycles can differ for each value of~$t$.
We emphasize that the completion of cycles is determined by the order of agents' operations,
and each agent can locally control the sequence of its own operations.
Hence, cycles can be optimized without making assumptions about
controlling the exact timing of agents' operations. 

\subsection{Optimizing over Constant~$c(t)$}
Suppose we wish to drive agents' iterates within a ball centered on the minimizer, $\hat{u}^{t}$, for all $t\in \mathcal{T}$.
In particular, for all~$t$, we would like agents' tracking error to be bounded 
by a given~$\rho > 0$ before the objective function changes from~$f(\cdot, t)$ to~$f(\cdot, t+1)$. 
This is first done for $c(t) \equiv c$ for all $t\in \mathcal{T}$. This corresponds to the 
case in which agents are able to do the same amount of communicating and 
computing for each objective function, as could happen if changes in the objective function were 
evenly spaced in time.
%\mh{Why do we do that?
%One explanation is: ``This corresponds to the case in which agents are, roughly speaking, able to do the same
%amount of communicating and computing for each objective function, as would happen if 
%changes in the objective function were evenly spaced in time.
%}
%\gb{done. Also could give a bound on how well the worst communicating agent needs to perform if~$c(t)=c_{min}$???}

\begin{theorem}[Finite Time Communication Cycle Requirement] \label{thm:thm3}
	%\mh{A theorem statement like this has to stand on its own.
	%So, define~$\rho$, define the error we want and say we want it to be less than~$\rho$.
	%Then say~$c(t) \equiv c$ for all~$t$, and, under these conditions, we attain
	%the desired bound if [the inequality below is true].
	%}	
	%\gb{done}
	Let Assumptions~\ref{as:c2}-\ref{as:tau} hold and
	let~$\rho > 0$. Suppose that N agents are executing Algorithm~\ref{alg:myAlg} with~$c(t) \equiv c$ for all~$t \in \mathcal{T}$. If
	\begin{equation} \label{eq:clim}
		c	\geq \frac{\ln\left(\frac{q_{max}^{\left(T+2\right)c} +\frac{\rho}{B}}{1+\frac{\rho}{B}}\right)}{\ln\left(q_{max}\right)},
	\end{equation}
	then $\left\Vert u^{i}(\eta_t)-\hat{u}^{t}\right\Vert _{2,\infty}\leq \rho$ for all~$t \in \mathcal{T}$. 
	%\gb{is it confusing to use~$c(T+2)$ in this result? Seems like it could refer to the number of cycles 2 time steps from now since that how we index~$c(t)$, however we do say that we are using a constant~$c$ instead of~$c(t)$ in this result.}
	%\mh{That can be confusing. Why not change the exponent to~$(T+2)c$?}	
	%\gb{done}
\end{theorem}
\emph{Proof:}    
	From Theorem~\ref{thm:thm2}, we have 	
	%\mh{Cite the theorem number instead of the equation number. So ``From Theorem 2...''}
	%\gb{done}
	\begin{equation}
		\left\Vert u^{i}(\eta_t)-\hat{u}^{t}\right\Vert _{2,\infty} \leq D_{0} \prod_{\theta=0}^{t} q_{\theta}^{c(\theta)} + \sum_{p=1}^{t} \sigma_{p-1} \prod_{r=p}^{t}q_{r}^{c(r)},
	\end{equation}
	where~$\eta_t = \sum_{\ell=0}^{t} \kappa_{\ell}$. For~$c(t) \equiv c$, this bound simplifies to
	%\mh{We need~$u^i(k)$ on the LHS, and we need to say which range of values of~$k$.}
	%\gb{took care of it by using~$\eta_{t}$}
	\begin{equation} \label{eq:cconst}
	    \left\Vert u^{i}(\eta_t)-\hat{u}^{t}\right\Vert _{2,\infty} \leq D_{0} \prod_{\theta=0}^{t} q_{\theta}^{c} + \sum_{p=1}^{t} \sigma_{p-1} \prod_{r=p}^{t}q_{r}^{c}.
	\end{equation}
	Equation~\eqref{eq:cconst} can be upper bounded via
	\begin{flalign}
	    \left\Vert u^{i}(\eta_t)-\hat{u}^{t}\right\Vert _{2,\infty} &\leq B \prod_{\theta=0}^{t} q_{max}^{c} + \sum_{p=1}^{t} B \prod_{r=p}^{t}q_{max}^{c} \\
	    %& \leq B \left( \prod_{\theta=0}^{t} q_{max}^{c} + \sum_{p=1}^{t} \prod_{r=p}^{t}q_{max}^{c} \right) \\
	    & \leq B \left( q_{max}^{(t+1)c} + q_{max}^{tc} + \ldots + q_{max}^{c} \right) \\
	    %& \leq Bq_{max}^{c} \left( q_{max}^{tc} + q_{max}^{(t-1)c} + \ldots + 1 \right) \\
	    %\left\Vert u^{i}(\eta_t)-\hat{u}^{t}\right\Vert _{2,\infty}
	    & \leq B q_{max}^{c} \sum_{\theta=0}^{t} q_{max}^{\theta c},
	\end{flalign}
	%\mh{where the first line [does something], the second line [does something]. Briefly explain here.
	%These can be short, i.e., ``the first line follows by factoring, the second line follows by re-arranging.'' No need to over-explain if it isn't needed.
	%}
	%\gb{done}
	where the first line follows from the definition of~$B$, 
	the second line expands the sum and products, and
	the third line combines like terms. 
	%the second line follows from factoring, 
	%the third line expands out the products and summation, 
	%the fourth line follows from factoring, 
	%and the last line is rewriting the previous line as the sum of a finite geometric series. 
	To bound the left-hand side by~$\rho$, it is sufficient to have
	\begin{equation} \label{eq:bsufficient}
	    B q_{max}^{c} \sum_{\theta=0}^{t} q_{max}^{\theta c} \leq \rho. 
	\end{equation}
	From the fact that~${\sum_{m=0}^{n-1} ar^{mc}=a\left(\frac{r^{cn}-1}{r^{c}-1}\right)}$ when~$r<1$, 
	the bound in Equation~\eqref{eq:bsufficient} holds if and only if 
	\begin{equation}
	    Bq_{max}^{c}\left(\frac{q_{max}^{(t+1)c}-1}{q_{max}^{c}-1}\right) \leq \rho. 
	\end{equation}	
	Setting~$t = T$ maximizes the left-hand side, and then 
	solving for~$c$ provides a bound that holds for all~$t \in \mathcal{T}$. 
	%\gb{I use~$T$ in the result for~$q_{max}^{c(T+2)}$, but I use~$t$ in my proof when I write~$q_{max}^{c(t+1)}$. Which should I use to be consistent? }
\hfill $\blacksquare$

%This lower bound on the required communication cycles can be regarded as a conservative bound
%since we consider the worst-case values of the problem parameters~$q_{max}$ and~$B$. 
Extending this bound to the asymptotic case gives the following. 

\begin{corollary}[Asymptotic Communication Cycle Requirement] \label{cor:cor2}
Let Assumptions~\ref{as:c2}-\ref{as:tau} hold and 
let~$\rho > 0$. Suppose that N agents are executing Algorithm~\ref{alg:myAlg} with~$c(t) \equiv c$ for all~$t \in \mathcal{T}$. 
If
	\begin{equation}
		c	\geq \frac{\ln\left(\frac{\frac{\rho}{B}}{1+\frac{\rho}{B}}\right)}{\ln\left(q_{max}\right)}
	\end{equation}	
	then $\left\Vert u^{i}(\eta_t)-\hat{u}^{t}\right\Vert _{2,\infty} \leq \rho$ for all~$t \in \mathbb{N}$.	
	%               \gb{should this be for all $t\in \N$, since its asymptotic?}
\end{corollary}

\emph{Proof:}
	In Theorem~\ref{thm:thm4}, use~$\lim_{T \to \infty} q_{max}^{c\left(T+2\right)} = 0$.
\hfill $\blacksquare$

%\mh{There's a subtle point here (we can discuss further, but I'll try to explain here).
%Here, agents' communications are informed by coupling in the cost~$f$. If I'm agent~$i$,
%then my gradient is~$\nabla_i f$, and my neighbors are whichever agents' states
%this depends on. So if~$\nabla_i f$ doesn't depend on, say,~$x_3$, then I never
%need to talk to agent~$3$. It's easy to show that this is symmetric ($i$ depends
%on~$j$ if and only if~$j$ depends on~$i$). What you say is completey true about
%some agents needing to send more, we should just be clear
%to the reader about why that is.
%}

%\begin{remark}
%\mh{What does this mean?}
%\gb{is this clear? We spoke about putting this result in before considering the structure of the network, do we still want it?}
%	\emph{How many messages does agent~$i$ send per communication cycle?}
%	\begin{equation}
%		m_i = \left|N_{i}\right|
%	\end{equation}
%	\mh{I see what this says, but what is this doing for us? Can we tie this back to the rest of the paper?}
%\end{remark}

\subsection{Optimizing over Time-Varying~$c(t)$}
Lastly, we optimize the number of communication cycles the network should complete for each objective function
to minimize dynamic regret subject to a constraint on the number of cycles agents can complete.
This can model, for example, the limits of finite battery life upon the operations an agent
can execute over the time horizon~$\mathcal{T}$. Although exact tracking error is not known
in general, we can minimize the bound on tracking error in 
Equation \eqref{eq:tracking_bound} because that bound is computable, and this is what
we will do.
%However, we are able to minimize the RHS of Equation~\eqref{eq:tracking_bound} 
%which will minimize~$\left\Vert u^{i}(k)-\hat{u}^{t}\right\Vert _{2,\infty}$. 
In particular, we solve the following optimization problem to calculate
the optimal number of cycles to complete for~$f(\cdot, t)$, denoted~$c^*(t)$. 

\begin{problem} \label{prob:cycle_opt}
%\mh{Give a complete problem statement: there are~$N$ agents, over the time horizon~$[T]$, etc.}
%\gb{done.}
For~$N$ agents running Algorithm~\ref{alg:myAlg} over a time horizon~$\mathcal{T}$,
    \begin{equation}
    	\begin{split}
    	\underset{\{c(t)\}_{t \in \mathcal{T}}}{\text{minimize}} & \ \sum_{t=0}^{T}q_{t}^{c(t)}D_{t} \\
    	\text{subject to} & \ \sum_{t=0}^{T}c(t)\leq K,
    	\end{split}
    \end{equation}
    %\mh{What is~$u$ doing here? Does this actually show up anywhere below?}
    %\gb{No it does not show up. I got rid of it.}
    where~$K \in \N$ is the total number of communication cycles agents can
    complete over the time horizon~$\mathcal{T}$. 
\end{problem}

Problem~\ref{prob:cycle_opt} is a constrained convex optimization problem whose
solution is derived next. 

\begin{theorem}[Optimal Number of Communication Cycles]\label{thm:thm4}
Let Assumptions~\ref{as:c2}-\ref{as:tau} hold, and define
%\mh{Again, state assumptions, conditions, etc. This theorem statement has to stand on its own. Say there are~$N$ agents,~$c(t)$ is the number of cycles to complete when minimizing~$f(\cdot, t)$, etc.}
    \begin{align}
    F_{t} &\coloneqq \sum_{\theta=t}^{T}\left(\prod_{i=0}^{t-1}q_{i}^{c^{*}(i)}\prod_{j=t+1}^{\theta}q_{j}^{c^{*}(j)}\right) \\
    E_{t} &\coloneqq \sum_{i=0}^{t-1}\sigma_{i}\sum_{j=1}^{T-t+1}\prod_{k=i+1}^{t-1}q_{k}\prod_{l=t+1}^{j+t-1}q_{l}. 
    \end{align}
    Then, for each~$t \in \mathcal{T}$, the optimal number of cycles for N agents executing Algorithm~\ref{alg:myAlg} to complete 
    when minimizing~$f(\cdot, t)$ 
    %(before it switches to~$f(\cdot, t+1)$) 
    is
    %between switches in objective functions from~$f(\cdot,t)$ to~$f(\cdot,t+1)$ is given by    
	\begin{equation}
		c^{*}(t) = \frac{\ln\left(-\frac{\mu^{*}}{\ln\left(q_{t}\right)\left[D_{0}F_{t}+E_{t}\right]}\right)}{\ln\left(q_{t}\right)}.
	\end{equation}
	%\mh{Don't use~$C_t$. We already have~$c(t)$ and~$C^i_j$. Just make it something else.}
	%\gb{done.}
	%\mh{What is~$\mu^*$?}
	%\gb{I don't have a closed-form solution for $\mu^{*}$. However, it can be solved for from the equations for~$c^{*}(t)$ along with the equation~$\sum_{t=0}^{T} c^{*}(t) = K$.}
\end{theorem}

%The solution of~$c^{*}(t)$ is explicit and depends on the values of~$c^{*}(t)$ at previous and future times. 
The value of~$\mu^*$ can be computed
along with~$c^*(t)$ by solving the system of equations
we derive at the end of the following proof. 

%\mh{It's fine to put brief remarks between a theorem statement and its proof. Here, say that~$c^*(t)$ is in terms of the~$c^*$ values at previous times.}
%\gb{done}

\emph{Proof:}	
	%Theorem \ref{thm:thm3} provides that
	%\begin{equation}\label{eq:proposed_sets2}
	%    \left\Vert u^{i}-\hat{u}^{t}\right\Vert _{2,\infty}\leq q_{t}^{c(t)}D_{t}.
	%\end{equation}
	%Therefore, minimizing the RHS of Equation \eqref{eq:proposed_sets2} will minimize the LHS.
    %
    Problem~\ref{prob:cycle_opt}'s objective function is twice continuously differentiable in each~$c(t)$.
    And, for each~$c(t)$, we find
    \begin{equation}
        \frac{\partial^2}{\partial c(t)^2} q_t^{c(t)} D_t = \big(\ln(q_t)\big)^2q_t^{c(t)} D_t > 0,
    \end{equation}
    and thus the objective function in Problem~\ref{prob:cycle_opt} is strictly convex. 
    %and satisfies the second-order condition for strict convexity,~${\nabla_{c(t)}^{2} \ q_{t}^{c(t)}D_{t} \succeq 0}$. Moreover, the sum operation of convex functions preserves convexity, therefore~${\sum_{t=0}^{T}q_{t}^{c(t)}D_{t}}$ is strictly convex wrt~$c(t)$. 
    %The constraint of Problem \ref{prob:cycle_opt} is affine in~$c(t)$, therefore the following Lagrangian function has a unique minimum,
    The Lagrangian associated with Problem~\ref{prob:cycle_opt} is
    \begin{equation} \label{eq:Lagrangian}
        L(c(t),\mu) = \sum_{t=1}^{T}q_{t}^{c(t)}D_{t} + \mu \left(\sum_{t=1}^{T} c(t) - K \right),
    \end{equation}
    where~$\mu \geq 0$ is the dual variable. 
    Because the constraint is affine, it satisfies the Linearity Constraint Qualification.
    Combined with the strict convexity of the cost, this implies that the Karush-Kuhn-Tucker (KKT) conditions
    are sufficient to find the unique minimizer for Problem~\ref{prob:cycle_opt}.
    
    To solve for the optimal primal-dual pair~$\left(c^{*}(t),\mu^{*}\right)$, we use the KKT first-order optimality conditions,~${\nabla_{c(t)}L(c(t)^{*},\mu^{*})=0}$  and~${\nabla_{\mu}L(c(t)^{*},\mu^{*})=0}$. 
    For a fixed~$t \in \mathcal{T}$, 
    computing the derivatives of~$L\big(c(t),\mu\big)$ with respect 
    to~$c(t)$ and~$\mu$ gives
    \begin{align}
        \frac{\partial L}{\partial c(t)} &= q_{t}^{c(t)} \ln(q_{t}) \left[D_{0}F_{t}+E_{t}\right] + \mu \\
        \frac{\partial L}{\partial \mu} &= \sum_{t=0}^{T} c(t) - K.
    \end{align}
    Using the KKT first-order optimality conditions, we set these derivatives equal to 0 for each~$t \in \mathcal{T}$ to find
    \begin{align}
        0 &= q_{0}^{c^{*}(0)} \ln(q_{0}) \left[D_{0}F_{0}+E_{0}\right] + \mu^{*} \\
        &\quad\quad\quad\quad\quad\quad\quad\vdots \\
        0 &= q_{T}^{c^{*}(T)} \ln(q_{T}) \left[D_{0}F_{T}+E_{T}\right] + \mu^{*} \\
        0 &= \sum_{t=0}^{T} c^{*}(t) - K.
    \end{align}
    %\mh{Check your sums over~$t$ throughout the paper. This one started at~$0$. I changed it to~$1$. Make sure that's the case everywhere.}
    %\gb{I've used~$t=0$ throughout my proofs in this paper. I think it would be easier to change the period~$[T] \coloneqq \{0,1,\ldots,T\}$ instead of changing my proofs?}
    %\mh{We could use~$\mathcal{T} = \{0, \ldots, T\}$.}
    %\gb{done}
    This is a system of~$T+2$ equations in~$T+2$ unknowns (namely,~$c^*(t)$ for~$t \in \mathcal{T}$ and~$\mu^* \in \mathbb{R}$). 
    Solving for~$c^{*}(t)$ and~$\mu^*$ completes the proof. 
    \hfill $\blacksquare$

	\section{Numerical Example}\label{sec:Simulation}
This section considers two classes of simulations. First, two time-varying quadratic programs
are considered to illustrate the results in Section~\ref{sec:results}.
Second, because the results of this paper are motivated in part by ongoing
work in feedback optimization (see the Introduction for discussion), the performance
of Algorithm~\ref{alg:myAlg} on a noisy feedback optimization problem is evaluated empirically.

\subsection{Time-Varying Quadratic Programs} \label{subsec:TV_QP}
We consider two examples with~$N = 15$ agents 
executing Algorithm~\ref{alg:myAlg} as their objective
switches~$10$ times. 
The objective function takes the form
\begin{equation}
    f(u,t) = \frac{1}{2} u^{T} Q(t) u + q(t)^{T} u,
\end{equation}
where~$0 \prec Q(t) = Q^T(t) \in \R ^{n \times n}$ 
and~$q(t) \in \R^{n}$ vary over $t\in \mathcal{T}$. Both~$Q(t)$ and~$q(t)$ are randomly generated in this example.
%\mh{How are you generating~$Q(t)$ and~$r(t)$? It's fine if they're random or just arbitrarily chosen.
%Make sure to tell the reader what's actually happening.}
%\gb{done}
At each tick of~$k$, agent~$i$ communicates with others with 
probability~$p_i(k) \in (0,1)$, for which
a new value is generated randomly from a uniform distribution at each~$k$. 
%$\mh{????} 
%\mh{What is that probability? Include it here.}
%\gb{done}
The objective function switches every 50 ticks of~$k$.
%\mh{Do these objective functions change every~$10,000$ ticks of~$k$? If so, say that.} \gb{done}

\begin{figure}
    \centering
    \includegraphics[width = 0.45 \textwidth]{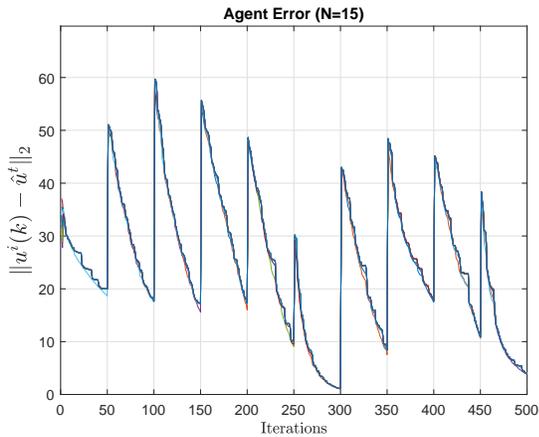}
    \caption{Plot of tracking error as a function of~$k$ for~$k$ from~$0$ to~$500$. 
    The value of~$\|u^i(k) - \hat{u}^t\|_2$ is shown for all~$i \in [N]$.
    %\mh{Regarding being a function of~$k$: it seems like it is.
    %Are you plotting distance to the optimizer at each~$k$? If so, the~$y$-axis label will
    %need to change. The same goes for the other plots.
    %}
    %\gb{done}
    %with~$N = 15$ agents and~$10$ switches in their objective function (i.e.,~$\mathcal{T} = \{0, \ldots, 10\}$). 
    %The periods of descent are when agents compute and communicate, leading to better 
    %tracking of their fixed point at that point in time, and the abrupt increases in 
    %tracking error are when their objective function changes.
    %While each change in the objective function causes
    %large increases in tracking error, agents are still able
    %to decrease tracking error before the next switch.
    }
    \label{fig:agent_error}
\end{figure}

\begin{figure}
    \centering
    \includegraphics[width = 0.45 \textwidth]{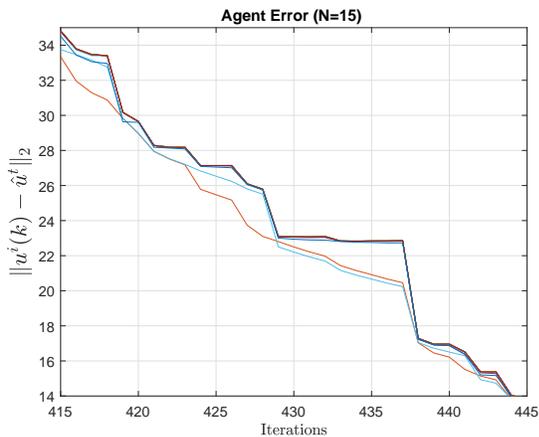}
    \caption{Zoomed-in plot of Figure~\ref{fig:agent_error} showing the tracking error as a function of~$k$ for
    all~$N = 15$ agents. 
    %\mh{Is this zooming in on Figure~1? If so, say that here.}
    %\gb{done}
    The periods of descent are when agents compute and communicate, 
    leading to better tracking of their fixed point at that time. 
    Flat areas are where agents either do not compute or communicate, 
    thereby completing no cycles and making no progress in reducing tracking error.
    }
    \label{fig:zoom_agents}
\end{figure}

Figure~\ref{fig:agent_error} shows how the distance between the current 
iterate~$u^{i}(k)$ and the unique minimizer~$\hat{u}^t$ changes for all~$i \in [N]$ 
and for each~$t \in \mathcal{T}$. It can be seen that the individual error of each agent decreases as more iterations
(and therefore more communication cycles) occur.
%\mh{It looks like Figure~1 sums over all agents' errors while Figure~2 is just for one choice of~$i$.
%Explain this more. I also think we could cut one of these. Either show Figure 1, which already
%includes every agent, or show just Figure 2 and say that it's representative of all agents' behaviors.
%}
Figure~\ref{fig:zoom_agents} 
is a zoomed-in plot of all agents' errors, and it 
shows the effects of total asynchrony on the convergence of Algorithm~\ref{alg:myAlg}.
In particular, the areas of descent in Figure~\ref{fig:zoom_agents} 
are due to agents computing and communicating, leading to completed
communication cycles and hence progress towards fixed points. 
Conversely, flat areas correspond to intervals in which cycles have not
yet been completed, thereby making no progress. 

After every~$50$ iterations there is a sudden increase in error because those are 
the times at which the objective function switches from~$f(\cdot,t)$ to~$f(\cdot,t+1)$. Thus, the 
latest copy of each agents' decision vector becomes the initial guess of the minimum for~$f(\cdot,t+1)$.
Because~$\hat{u}^{t+1} \neq \hat{u}^t$, this initial guess is typically far away, leading to
an abrupt increase in error at switching times. 


\subsection{Feedback Optimization with Measurement Noise}
Secondly, we consider an example from the study of \emph{feedback optimization}~\cite{colombino2019,bernstein2019,menta2018}
%\gb{I added citations here for feedback optimization}
in which the cost is on the deviation of the output of a dynamical system from a reference trajectory.
Consider the objective function
\begin{equation}
    f(u,t) = \frac{1}{2} u^{T}Q(t)u + \frac{1}{2} \left(y(u) - r(t) \right)^{2},
\end{equation}
%\mh{We need a notation change here. Above~$r(t)$ was just a vector in the cost, but now it's a reference signal.}
%\gb{I changed the first example's~$r(t) \rightarrow p(t)$ and kept the reference as~$r(t)$}
where~$Q(t) \in \R^{n \times n}$, $y:\R^{n} \rightarrow \R$ is the input-output map of the system, 
and~$r(t) \in \R$ is the time-varying reference signal that~$y$ must track. 
In feedback optimization we only have knowledge of an output measurement affected by noise. Therefore, the objective function that we can compute is
\begin{equation}
    f(u,t) = \frac{1}{2} u^{T}Q(t)u + \frac{1}{2} \left( \hat{y}(u) - r(t) \right)^{2},
\end{equation}
where $\hat{y}(u)=y(u)+n(t)$ is the noisy measurement of the input-output map and~$n(t)\sim\mathcal{N}(0,{1000^2})$ is i.i.d.. 
%\gb{I said what kind of noise I used.}

\begin{figure}
    \centering
    \includegraphics[width = 0.45 \textwidth]{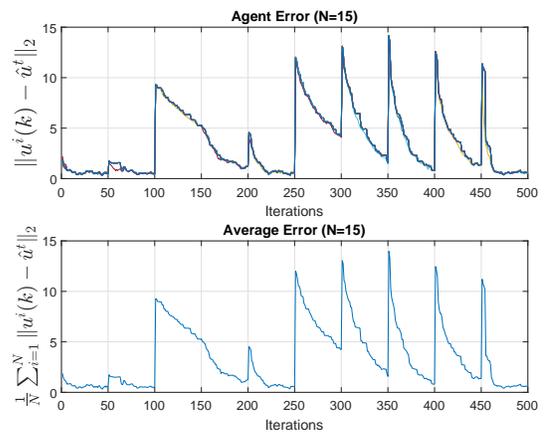}
    \caption{Plot of tracking error as a function of~$k$ 
    for all~$N = 15$ agents in the presence of noisy output measurements leading to a noisy gradient calculation. The sudden increases in error are due to the switching of objective functions similar to the previous example. Effects of asynchrony and noisy measurements can be observed; however, agents are still able to track their fixed point for all~$t$ within a
    relatively small ball centered on them.
    }
    \label{fig:feedback_opt}
\end{figure}

Figure~\ref{fig:feedback_opt} shows that, even with measurements corrupted by noise, Algorithm~\ref{alg:myAlg} 
can still track time-varying fixed points to within a ball centered on them. For example, just before the~$8^{th}$ switch, agents' average tracking
error is 1.192, which indicates close tracking even under noise.
%\mh{Complete this sentence: ``For example, just before the~$8^{th}$ switch, agents' total tracking
%error is \#\#\#\#, which indicates close tracking even under noise.''}
Similar to the previous example, the switching of objective functions causes a noticeable error that is quickly 
mitigated by subsequent communication cycles. 
This example goes beyond the analytical tracking guarantees provided above, and further
characterization of Algorithm~\ref{alg:myAlg} with measurement noise is a direction
for future research.

	\section{Conclusions} \label{sec:conclusions}
In this paper, we have proposed a totally asynchronous multi-agent algorithm for 
tracking time-varying fixed points that are solutions to time-varying optimization
problems. We further optimized agents' computations and communications to minimize dynamic regret
when they are able to perform only a limited number of each
type of operation. Next steps in this work include integrating the dynamics of 
agents into the update law for feedback optimization and analytically bounding tracking
error when measurement noise is included in a feedback optimization setup.

	%%%%%%%%%%%%%%%%%%%%%%%%%%%%%%%%%%%%%%%%%%%%%%%%%%%%%%%%%%%%%%%%%%%%%%%%%%%%%%%%	
	\bibliographystyle{IEEEtran}{}
	\bibliography{sources}
	\appendix

\subsection{Proof of Proposition \ref{prop:TI_convergence}}
\label{app:TI_convergence}
%\gb{do we need this appendix if we are citing a CORE Lab paper for the proof?}
%\mh{Leave this here for now. Given that we're over the page limit, we'll upload a version
%of this (with this appendix included) to arXiv. Then when we submit to the conference,
%we'll cut this appendix from the submitted version and cite it on
%arXiv.}
%\gb{okey-dokey}
%\gb{I replaced all the~$D_{t}$ with~$\left\Vert u^{i}(\eta_{t-1}) - \hat{u}^t \right\Vert_{2,\infty}$}                                                    
We prove that each of the four properties of Definition~\ref{def:sets} is satisfied.
    First, consider the sets
	\begin{equation}
	    U_{i}(s) \!=\! \left\{\!y_{i}\in\mathbb{R}^{n_{i}} \!:\! \left\Vert y_{i}-\hat{u}_{i}^{t}\right\Vert _{2}\leq q_{t}^{s}\max_{i \in [N]}\left\Vert u^{i}(\eta_{t-1}) - \hat{u}^t \right\Vert_{2,\infty}\!\right\}.
	\end{equation}
	Then, for~$y = (y_1, \ldots, y_N)^T \in U_1(s) \times \cdots U_N(s)$, we have~$N$
	inequalities
    \begin{align}
    \left\Vert y_{1}-\hat{u}_{1}^{t}\right\Vert _{2} &\leq q_{t}^{s}\max_{i \in [N]}\left\Vert u^{i}(\eta_{t-1}) - \hat{u}^t \right\Vert_{2,\infty} \\
    \quad\quad\quad\quad\vdots \\
    \left\Vert y_{N}-\hat{u}_{N}^{t}\right\Vert _{2} & \leq q_{t}^{s}\max_{i \in [N]}\left\Vert u^{i}(\eta_{t-1}) - \hat{u}^t \right\Vert_{2,\infty}, 
    \end{align}
	from which it follows that
	\begin{equation}
	    \left\Vert y-\hat{u}^{t}\right\Vert _{2,\infty} \coloneqq \max_{i \in [N]} \left\Vert y_{i}-\hat{u}_{i}^{t}\right\Vert _{2}\leq q_{t}^{s} \max_{i \in [N]} \left\Vert u^{i}(\eta_{t-1}) - \hat{u}^t \right\Vert_{2,\infty}.
	\end{equation}
	Thus~$y\in U_{1}(s)\times \cdots \times U_{N}(s)$ implies~$y\in U(s)$ and
	\begin{equation} \label{eq:containment1}
	    U_{1}(s)\times U_{2}(s)\times\cdots\times U_{N}(s) \subseteq U(s).
	\end{equation}
	Now, consider~$y\in U(s)$. Then,
	    \begin{align}
	        \left\Vert y-\hat{u}^t\right\Vert _{2,\infty}\leq q_{t}^{s}\max_{i \in [N]} \left\Vert u^{i}(\eta_{t-1}) - \hat{u}^t \right\Vert_{2,\infty} 
	   \end{align}
	   and
	   \begin{align}
	        \max_{i \in [N]} \left\Vert y_{i}-\hat{u}_{i}^{t}\right\Vert _{2}\leq q_{t}^{s}\max_{i \in [N]}\left\Vert u^{i}(\eta_{t-1}) - \hat{u}^t \right\Vert_{2,\infty}.
	    \end{align}
	    This implies that, for all~$i$, we have
	    %\begin{equation}
	        ${\left\Vert y_{i}-\hat{u}_{i}^{t}\right\Vert _{2} \leq q_{t}^{s} \max\limits_{i \in [N]} \left\Vert u^{i}(\eta_{t-1}) - \hat{u}^t \right\Vert_{2,\infty}}$.
	    %\end{equation}
	    Then~$y\in U(s)$ implies that~$y_{i}\in U_{i}(s)$ for all~$i$, and thus $U(s) \subseteq U_{1}(s)\times U_{2}(s)\times\cdots\times U_{N}(s)$.
	    This, combined with Equation~\eqref{eq:containment1}, implies that Definition~\ref{def:sets}.1 holds. 
	 
	 For Definition~\ref{def:sets}.2, consider~$y \in U(s+1)$. Then
	    \begin{equation}
	        \left\Vert y-\hat{u}^{t}\right\Vert_{2,\infty} \leq q_{t}^{s+1} \max_{i \in [N]} \left\Vert u^{i}(\eta_{t-1}) - \hat{u}^t \right\Vert_{2,\infty}.
	    \end{equation}
	    Since~$q_{t}\in (0,1)$, we have
	    \begin{multline}
	        \left\Vert y-\hat{u}^{t}\right\Vert _{2,\infty} \leq q_{t}^{s+1} \max_{i \in [N]} \left\Vert u^{i}(\eta_{t-1}) - \hat{u}^t \right\Vert_{2,\infty} \\ < q_{t}^{s} \max_{i \in [N]} \left\Vert u^{i}(\eta_{t-1}) - \hat{u}^t \right\Vert_{2,\infty},
	    \end{multline}
	    which implies~$y\in U(s)$. Then~$U(s+1)\subseteq U(s)$.
	  
	For Definition~\ref{def:sets}.3, we note that, because~$q_t \in (0, 1)$, we have~$\lim\limits_{s \to \infty} q_t^s = 0$.
	It immediately follows from the definition of~$U(s)$ that
	%\begin{equation}
	$\lim_{s \to \infty} U(s) = \{\hat{u}^t\}$. 
	%\end{equation}

    Towards satisfying Definition~\ref{def:sets}.4, we first observe that~\cite[Lemma 1]{ubl2020} implies
	    \begin{multline} \label{eq:2infbound}
	        \left\Vert I_{n}-\gamma_{t}H(u,t)\right\Vert _{2,\infty} \\ \leq \max_{i \in [N]} \left[\left\Vert I_{n_{i}}-\gamma_{t}H_{ii}(u,t)\right\Vert_{2}
	        +\gamma_{t}\underset{j \neq i}{\sum_{j=1}^{n}}\left\Vert H_{ij}(u,t)\right\Vert _{2}\right].
	    \end{multline}

    %Let~$H(u,t)$ satisfy Assumption \ref{as:diag}. \mh{This should be a hypothesis in the proposition statement.}
	    %We want to design a~$\gamma_{t}$ that makes the following inequality true for all~${u\in U}$ and for each~$t\in\left[ T \right]$
	    %\begin{equation}
	    %    \left\Vert I_{n}-\gamma_{t} H(u,t) \right\Vert _{2,\infty} \leq 1-\gamma_{t} \beta_{t}.
	    %\end{equation}
	 %   First we note \cite[Lemma 1]{ubl2020} implies
	 %   \begin{flalign}
	 %       &\left\Vert I_{n}-\gamma_{t}H(u,t)\right\Vert _{2,\infty} \\ &\leq \max_{i \in [N]} \left[\left\Vert I_{n_{i}}-\gamma_{t}H_{ii}(u,t)\right\Vert _{2}+\gamma_{t}\underset{j \neq i}{\sum_{j=1}^{n}}\left\Vert H_{ij}(u,t)\right\Vert _{2}\right].
	 %   \end{flalign}
	    %Therefore, it is sufficient to upper bound the RHS of this equation by~${1-\gamma_{t} \beta_{t}}$ which we show next. First note the following
	    Now consider the~$i^{th}$ diagonal block of~$I_n - \gamma_tH(u, t)$, which satisfies %has size~$n_i \times n_i$ and satisfies
	    %\begin{equation}
	        $\left[I_{n}-\gamma_{t}H(u,t)\right]_{ii}=I_{n_{i}}-\gamma_{t}H_{ii}(u,t)$.
	    %\end{equation}
	    By Assumption \ref{as:diag},~$H_{ii}(t)$ is symmetric. Then
	    \begin{multline}
	        \|I_{n_i}-\gamma_{t}H_{ii}(u,t)\|_{2}\\
	        =\max\Big\{\!\big|\lambda_{min}[I_{n_{i}} \!-\! \gamma_{t}H_{ii}(u,t)]\big|,\big|\lambda_{max}[I_{n_{i}} \!-\! \gamma_{t}H_{ii}(u,t)]\big|\!\Big\} \\
	        =\max\Big\{\!\big|1-\gamma_{t}\lambda_{min}[H_{ii}(u,t)]\big|,\big|1-\gamma_{t}\lambda_{max}[H_{ii}(u,t)]\big|\Big\}.
	    \end{multline}
	    %We would like to design~$\gamma_{t}$ such that the following remains true for all~$u\in U$ and for each~$t\in [T]$
	    %\begin{align}
	    %    \gamma_{t}&<\frac{1}{\underset{u\in U}{max}\ \underset{i\in\left[N\right]}{max}\lambda_{max}\left(H_{ii}(u,t)\right)}\\&=\frac{1}{\underset{u\in U}{max}\ \underset{i\in\left[N\right]}{max}\left\Vert H_{ii}(u,t)\right\Vert _{2}}\\&\leq\frac{1}{\underset{u\in U}{max}\ \underset{i\in\left[N\right]}{max}\underset{j \neq i}{\sum_{j=1}^{n}}\left\Vert H_{ij}(u,t)\right\Vert _{2}}
	    %\end{align}
	    The upper bound on~$\gamma_t$ gives ${1-\gamma_{t} \lambda_{min}\left[ H_{ii}(u,t) \right] \geq 1-\gamma_{t}\lambda_{max} \left[ H_{ii}(u,t) \right]}$. 
	    Then
	    \begin{equation}
	        \left\Vert I_{n_{i}}-\gamma_{t}H_{ii}(u,t)\right\Vert _{2} =1-\gamma_{t}\lambda_{min}[H_{ii}(u,t)].
	    \end{equation}
	    Then we have 
	    \begin{flalign}
	        &\underset{i\in\left[N\right]}{\max} \,\,\left\Vert I_{n_{i}}-\gamma_{t}H_{ii}(u,t)\right\Vert _{2}+\gamma_{t}\underset{j \neq i}{\sum_{j=1}^{n}}\left\Vert H_{ij}(u,t)\right\Vert_{2} \\
	        %&=\underset{i\in\left[N\right]}{max}\left[1-\gamma_{t}\lambda_{min}\left(H_{ii}(u,t)\right)+\gamma_{t}\underset{j \neq i}{\sum_{j=1}^{n}}\left\Vert H_{ij}(u,t)\right\Vert _{2}\right]\\
	        &=\!\underset{i\in\left[N\right]}{\max} \,\,1 \!-\! \gamma_{t} \Big(\lambda_{min}[H_{ii}(u,t)] \!-\! \underset{j \neq i}{\sum_{j=1}^{n}}\left\Vert H_{ij}(u,t)\right\Vert _{2}\!\Big) \\
	        &\leq1-\gamma_{t}\beta_{t}. 
	    \end{flalign}
	    Using Equation~\eqref{eq:2infbound},
	    %\begin{equation} \label{eq:setsmain}
	        $\left\Vert I_{n}-\gamma_{t}H(u,t)\right\Vert_{2,\infty} \leq1-\gamma_{t}\beta_{t}$. 
	    %\end{equation}
	    
	    For Definition~\ref{def:sets}.4, consider
	    %\begin{equation}
	        $\hat{u}^{t} =\Pi_{U}\left[\hat{u}^{t}-\gamma_{t}\nabla_{u}f\left(\hat{u}^{t}, t\right)\right]$,
	    %\end{equation}
	    i.e.,~$\hat{u}^t$
	    is a fixed point of the projected gradient descent mapping. Using the non-expansive property of the projection operator~$\Pi_{U}\left[\cdot\right]$, we find
	    \begin{multline}
	        \left\Vert z_{i}-\hat{u}_{i}^{t}\right\Vert_{2} \\
	        =\left\Vert \Pi_{U_{i}}\left[y_{i}-\gamma_{t}\nabla_{u_{i}}f(y, t)\right]-\Pi_{U_{i}}\left[\hat{u}_{i}^{t}-\gamma_{t}\nabla_{u_{i}}f\left(\hat{u}^{t}, t\right)\right]\right\Vert _{2}\\
	        \leq\left\Vert y_{i}-\gamma_{t}\nabla_{u_{i}}f(y,t)-\hat{u}_{i}^{t}+\gamma_{t}\nabla_{u_{i}}f\left(\hat{u}^{t}, t\right)\right\Vert_{2}.
	    \end{multline}
	    Using the Mean Value Theorem~\cite{flett1958}, we have  
	    \begin{equation}
	        \nabla f(y, t) - \nabla f(\hat{u}^{t}, t) %&= \int_{0}^{1} H\big(\hat{u}^{t} + \tau (y-\hat{u}^t), t\big) d\tau \cdot (y-\hat{u}^t) \\ 
	        = M(y,t) (y-\hat{u}^t),
	    \end{equation}
	    where~$M(y,t) \coloneqq \int_{0}^{1} H(\hat{u}^{t} + \tau (y-\hat{u}^t),t) d\tau$, which then gives
        \begin{equation}
            \left\Vert I_{n}-\gamma_{t}M(y,t)\right\Vert _{2,\infty} \leq1-\gamma_{t}\beta_{t}.
        \end{equation}
        Then we have
        \begin{flalign}
            &\left\Vert y_{i}-\gamma_{t}\nabla_{u_{i}}f(y,t)-\hat{u}_{i}^{t}+\gamma_{t}\nabla_{u_{i}}f\left(\hat{u}^{t},t\right)\right\Vert _{2}\\
            &\leq\underset{i\in\left[N\right]}{\max}\left\Vert y_{i}-\gamma_{t}\nabla_{u_{i}}f(y,t)-\hat{u}_{i}^{t}+\gamma_{t}\nabla_{u_{i}}f\left(\hat{u}^{t},t\right)\right\Vert_{2}\\
            &=\left\Vert y-\gamma_{t}\nabla_{u}f(y,t)-\hat{u}^{t}+\gamma_{t}\nabla_{u}f\left(\hat{u}^{t},t\right)\right\Vert_{2,\infty}\\
            &=\left\Vert y-\hat{u}^{t}-\gamma_{t}\left(\nabla_{u}f(y,t)-\nabla_{u}f\left(\hat{u}^{t},t\right)\right)\right\Vert_{2,\infty}\\
            &=\left\Vert y-\hat{u}^{t}-\gamma_{t}M(y,t)\left(y-\hat{u}^{t}\right)\right\Vert_{2,\infty}\\
            &\leq\left\Vert I_{n}-\gamma_{t}M(y,t)\right\Vert _{2,\infty}\left\Vert y-\hat{u}^{t}\right\Vert_{2,\infty}\\
            &\leq\left(1-\gamma_{t}\beta_{t}\right)\left\Vert y-\hat{u}^{t}\right\Vert_{2,\infty}.
        \end{flalign}
        Therefore, because the bounds on~$\gamma_t$ imply~$1 - \gamma_t\beta_t \in (0, 1)$, 
        $y\in U(s)$ implies that~${z_{i}\in U_{i}(s+1)}$.

\end{document}